\documentclass[12pt,reqno]{amsart}
\usepackage{amsmath,amsthm,amsbsy,amssymb, longtable,hhline}%
\topmargin by-\headheight \advance
\topmargin by-\headsep
\evensidemargin
\oddsidemargin \marginparwidth 0.5in
\usepackage{amsmath, amssymb, graphicx, array, verbatim, %
caption, ccaption, %
epsfig, xspace}

\newtheorem{theorem}{Theorem}%
\newtheorem{lemma}{Lemma}%
\newtheorem{corollary}{Corollary}%
\usepackage{footnpag} %

\newcommand{\Area}{\mathrm{Area}}
\newcommand{\length}{\mathrm{length}}

\newcolumntype{C}{>{$}c<{$}}
\newcolumntype{L}{>{$}l<{$}}
\newcolumntype{R}{>{$}r<{$}}
\def\A{\mathcal A} 

\newcommand{\N}{{\mathbf{N}}}
\newcommand{\T}{\mathcal{T}}

\newcommand{\I}{\mathcal{I}}

\newcommand{\kk}{{\mathbf k}}

\newcommand{\ff}{{\mathbf f}}
\newcommand{\qq}{{\mathbf q}}

\def\RR{\mathbb{R}}
\def\NN{\mathbb{N}}
\def\d{\mathfrak d}

\def\cc{\mathbf c}

\newcommand{\F}{{\mathfrak{F}}}
\newcommand{\cF}{c\!\mathfrak{F}}
\newcommand{\cD}{c\!\mathfrak{D}}

\newcommand{\G}{{\mathfrak{G}}}
\newcommand{\K}{{\mathcal{K}}}
\newcommand{\E}{{\mathcal{E}}}

\newcommand{\FQ}{{\mathfrak{F}_{_Q}}}

\newcommand{\FQo}{{\mathfrak{F}_{_{Q,\mathrm{odd}}}}}

\newcommand{\norm}[1]{\lVert#1\rVert}

\def\smod{\!\!\!\!\! \pmod}

\begin{document}

\author[C. Cobeli and A. Zaharescu]
{Cristian Cobeli and Alexandru Zaharescu}\footnotetext{CC is partially supported by the CERES Programme of the Romanian Ministry of Education and Research, contract 4-147/2004.}

\address{
CC:
Institute of Mathematics of the Romanian Academy,
P.O. Box \mbox{1-764}, Bucharest 70700,
Romania.}
\email{cristian.cobeli@imar.ro}

\address{
AZ:
Institute of Mathematics of the Romanian Academy,
P.O. Box \mbox{1-764}, Bucharest 70700,
Romania.}

\curraddr{
AZ:
Department of Mathematics,
University of Illinois at Urbana-Champaign,
Altgeld Hall, 1409 W. Green Street,
Urbana, IL, 61801, USA.}
\email{zaharesc@math.uiuc.edu}

\title[The distribution of rationals in residue classes]
{The distribution of rationals in residue classes}

\subjclass[2000]{Primary 11B57}
\thanks{Key Words and Phrases: Farey fractions, parity problem}

\begin{abstract}
Our purpose is to give an account of the $r$-tuple problem on the
increasing sequence of reduced fractions having denominators bounded
by a certain size and belonging to a fixed real interval. We show that
when the size grows to infinity, the proportion of the $r$-tuples of
consecutive denominators with components in certain apriori fixed
arithmetic progressions with the same ratio approaches a limit, which is
independent on the interval. The limit is given explicitly and it is
completely described in a few particular instances. 
\end{abstract}

\maketitle

\section{Introduction}\label{Section1}

Let $Q$ be a positive integer and let $\I$ be an interval of real
numbers. We denote by $\F_Q^\I$ the sequence of reduced fractions from
$\I$, whose denominators are positive and $\le Q$. The elements of the sequence
are assumed to be arranged in ascending order. 
Since the denominators of these fractions are periodic with respect to an unit
interval, and they also determine uniquely the numerators, one usually
focuses on $\FQ$, the sequences corresponding to $\I=[0,1]$. This is known
as the Farey sequence of order $Q$.

Questions concerned with Farey sequences have a long history.
In some problems, such as those related to the 
connection between Farey fractions and Dirichlet $L-$functions, 
one is lead to consider subsequences of Farey fractions defined 
by congruence constraints. Knowledge of the distribution
of subsets of Farey fractions with congruence constraints
would also be useful in the study of the
periodic two-dimensional Lorentz gas. This is a billiard
system on the two-dimensional torus with one or more
circular regions (scatterers) removed 
(see Sina\u {\i}~\cite{Si}, Bunimovich~\cite{Bu}, Chernov~\cite{Ch},
Boca and Zaharescu~\cite{BZ}). 
Such systems were introduced in 1905 by 
Lorentz~\cite{L} to describe the dynamics of
electrons in metals. A problem raised by Sina\u {\i}
on the distribution of the free path length for this
billiard system, when small scatterers are placed at integer
points and the trajectory of the particle
starts at the origin, was solved in Boca et all~\cite{BGZ1}, \cite{BGZ2},
using techniques developed in ~\cite{ABCZ}, \cite{BCZ1}, 
\cite{BCZ2} to study the local spacing distribution of Farey 
sequences. 

The more general case when the trajectory starts 
at a given point with rational coordinates is intrinsically 
connected with the problem of the distribution of
Farey fractions satisfying congruence constraints. 
For example, the case when the trajectory starts from the center 
$(1/2,1/2)$ of the unit square
is related to the distribution of Farey fractions
with odd numerators and denominators. 

Some questions on the distribution of Farey fractions with
odd, and respectively even denominators have been investigated 
in~\cite{BCZ3}, \cite{Integers}, \cite{Evens}, and~\cite{H}. Most
recently, the authors~\cite{CZcd} have proved the existence of a density
function of neighbor denominators that are in an arbitrary arithmetic
progression. 

Let $d\ge 2$ and $0\le c \le d-1$ be integers.
We denote by $\FQ(c,d)$ the set of Farey fractions of order $Q$ with
denominators $\equiv c \pmod d$. We assume that the elements
of $\FQ(c,d)$ are arranged increasingly.
In various problems on the distribution of Farey fractions
with denominators in a certain arithmetic progression it would be very
useful to understand how the set $\FQ(c,d)$ sits inside $\FQ$. 
With this in mind, in the present paper we investigate the
distribution modulo $\d$ of the components of tuples of consecutive Farey 
fractions. To make things precise, let $\cc=(c_1,\dots,c_r)\in
[0,\d-1]^r\cap \NN^r$, for some integer $r\geq1$.
We say that a tuple $\qq=(q_1,\dots,q_r)$ of consecutive denominators
in $\FQ$ has \emph{parity} $\cc$ if $q_j\equiv c_j \pmod \d$,
for $1\le j\le r$. In this case, we write shortly $\qq\equiv \cc\pmod
\d$. We also say that $\cc$ is the \emph{parity}  
of the tuple of fractions $\ff=(a_1/q_1,\dots,a_r/q_r)\in\F_Q^r$.
The question, we address here, is whether it is true that
$\rho_Q(\cc,\d)$, the proportion of $r-$tuples of consecutive Farey
fractions  $\ff\in\F_Q^r$ of parity $\cc$, has a limit $\rho(\cc,\d)$, as
$Q\rightarrow\infty$. And, if so, can one provide an  explicit formula
for $\rho(\cc,\d)$\,? More generally, the same questions can  be asked for the
subset of Farey fractions which belong to a given real interval.

These problems have positive answers, indeed. Thus, for any interval
of positive length $\I$, any $r\geq1$ and any choice of parities
$\cc$, the proportion of $r-$tuples of consecutive 
fractions in $\F_Q^{\I}$ of parity  $\cc$ approaches 
a limit as $Q\rightarrow\infty$. Moreover, this limit depends
on $r$, $\cc$ and $\d$ only, and is independent of the choice of the
interval $\I$. Roughly, this says that the probability that $r$
consecutive fractions have parity $\cc$ is independent of the
position of these fractions in $\RR$. 

For instance, when $r=1$, one finds that the probability that a randomly
chosen fraction is even\footnote{We call a fraction \emph{odd}
or \emph{even} according to whether its denominator is odd or even.}
equals $1/3$. 
In other words, there are asymptotically 
twice as many odd fractions as even fractions in $\FQ$.
For $r\geq2$, the parities of $r$ consecutive fractions
are not independent of one another. Thus, for instance, the
probability that two neighbor fractions are even is zero.
This follows by a classical fundamental property, which says that any
consecutive Farey fractions $a'/q', a''/q''$ satisfy
	\begin{equation}\label{eqf}
		a''q'-a'q''=1.
	\end{equation}
This implies that $q',q''$ are relatively prime, so not both of them are even. 

In Section~\ref{SecPere} we present some properties of $\FQ$ that will
be used in the proofs. The next sections are devoted to give precise
statements of the main results and their proofs. We conclude with a few
examples for some special arithmetic progressions.

\section{Facts about $\FQ$}\label{SecPere}
Here we state a few  fact about Farey sequences that are needed in the
sequel. For the proofs and in depth details we refer to~\cite{BCZ1},
\cite{BCZ2}, \cite{Tessellations}, \cite{Hall70} and \cite{HT}.

For any $r\ge 1$, let $\cF_Q^r$ and $\cD_Q^r$ be the set of $r$-tuples
of consecutive fractions in $\FQ$, respectively denominators of
fractions in $\FQ$.  

We begin with a geometric correspondent of $\FQ$. Let $\T$ be the
triangle with vertices $(1,0)$, $(1,1)$, $(0,1)$ and its scaled
transform, $\T_Q:=Q\T$. Any lattice point with coprime coordinates
$(q',q'')\in Q\T$ determines uniquely a pair of consecutive
fractions $(a'/q',a''/q'')\in\cF_Q$, and conversely, if
$(a'/q',a''/q'')\in\cF_Q$ then $\gcd(q',q'')=1$, $q'+q''>Q$, so
$(q',q'')\in Q\T$. 

Next, let us notice that one way to read \eqref{eqf} yields
$a' = -\overline{q''}\pmod {q'}$ and 
$a' = \overline{q'}\pmod {q'}$. 
(Here the representative of $\overline x$, the inverse of $x\pmod d$,
is taken in the interval $[0,d-1]$.)  
This is the reason and the main idea which supports the fact many
properties regarding the distribution of $\FQ$ are preserved on
$\FQ\cap \I$, for any $\I\subseteq [0,1]$, and then by periodicity for
all larger intervals $\I$.

Another noteworthy property says that, starting with a pairs
$(q',q'')\in\cD^2_Q$, then one may calculate by a recursive method all
the fractions that precede or succeed $a'/q'$ and $a''/q''$. To see
this, we first remark that if $a'/q'$ and $a''/q''$ are not neighbor
in $\FQ$, then \eqref{eqf} no longer holds true. However, there is a
good replacement. Indeed, if $(a'/q',a''/q'',a'''/q''')\in\cF^3_Q$,
then the \emph{median} fraction is linked to the extremes through the
integer 
	\begin{equation}\label{eqmediumk}
		k=\frac{a'+a'''}{a''}=\frac{q'+q'''}{q''}
			=\Big[\frac{q'+Q}{q''}\Big]\,.
	\end{equation}
Suppose now that $r\ge 3$ and $\qq=(q_1,\dots,q_r)\in\cD^r_Q$, and
consider the positive integers defined by
	\begin{equation*}
		k_j=\frac{q_{j}+q_{j+2}}{q_{j+1}}
			=\Big[\frac{q_{j}+Q}{q_{j+1}}\Big]\quad
		\text{ for }j=1,\dots, r-2\,.
	\end{equation*}
Then we may say that $(q_1,q_2)$ \emph{generates} $\qq$ and
$\kk=(k_1,\dots,k_{r-2})$, too. In order to emphasize this, we write
	\begin{equation*}
		\begin{split}
		\qq(q_1,q_2)&=\qq^r_Q(q_1,q_2)=(q_1,\dots,q_r),\\
		\kk(q_1,q_2)&=\kk^{r-2}_Q(q_1,q_2)=(k_1,\dots,k_{r-2}),
		\end{split}
	\end{equation*}
dropping the sub- and superscript when they are clear from the
context. 

For a given positive integer $s$, let $\A^s$ be the set of all
\emph{admissible} $\kk\in(\NN^*)^s$, which we define to be those
vectors $\kk$ generated by some $(q',q'')\in\T_Q$, for some $Q\ge 1$ with
$\gcd(q',q'')=1$. We remark that $\A^s\subsetneq (\NN^*)^s$, and the
components of vectors from $\A^s$ should relate to one another in
specific manners, such as: 
a neighbor of $1$ can be any integer $\ge 2$,
neighbors of $2$ can be only $1,2,3$ or $4$,
neighbors of $3$ can be only $1$ or $2$, 
neighbors of $4$ can be only $1$ or $2$,
and a neighbor of any $k\ge 5$ can be only $1$.
Though, these conditions do not characterize completely the elements
of $\A^s$ for $s\ge 3$, and the complexity grows with $r$.

Although $\qq(q_1,q_2)$ is uniquely generated by its first two
components $q_1,q_2$, the vectors $\kk\in\A^s$ have each infinitely
many generators. The set of generators of a $\kk\in\A^s$ can be nicely
described and it will play a special role in what follows. We do this
by introducing its natural continuum envelope downscaled by a factor
of $Q$. Let us see the precise definitions. For any $(x,y)\in \RR^2$,
let $\big\{L_j(x,y)\big\}_{j\geq -1}$ be the sequence given by
$L_{-1}(x,y)=x$, $L_0(x,y)=y$ and recursively, for $j\geq 1$,
    \begin{equation*}
        L_j(x,y)=\left[ \frac{1+L_{i-2}(x,y)}{L_{i-1}(x,y)} \right]
        L_{i-1}(x,y)-L_{i-2}(x,y)\,.
    \end{equation*}
Then, we put
    \begin{equation*}
      \kk :\T \to (\NN^*)^{s} ,\qquad
       \kk(x,y)=\big( k_1(x,y),\dots ,k_{s} (x,y)\big),
    \end{equation*}
where, for $1\le j\le s$,
      \begin{equation}\label{eqk}
         k_j(x,y)=\left[ \frac{1+L_{j-2}(x,y)}{L_{j-1}(x,y)} \right]. 
      \end{equation}
Now, for any $\kk \in (\NN^*)^{s}$, we define
       \begin{equation*}
	   \T[\kk]=\T^s[\kk] :=\big\{ (x,y)\in \T \colon\ 
              \kk(x,y)=\kk \big\} \,,
       \end{equation*}
the domain on which the map $\kk(x,y)$ is locally constant.
The definition may be extended for an empty $\kk$ (with no
components) by putting $\T^0[\cdot]=\T$, the Farey triangle.
It turns out that $\T[\kk]$ is always a convex polygon.
Also, for any fixed $s\ge 0$, the set of all polygons
$\T[\kk]$, with $\kk\in\A^s$, form partition of $\T$,
that is, $\T =\bigcup\limits_{\kk\in \A^s} \T[\kk]$
and $\T[\kk] \cap \T[\kk'] =\emptyset$, whenever $\kk,\kk'\in\A^s$,
$\kk \neq \kk'$. 

The set we were looking for is the polygon $Q\T[\kk]$, which contains all
lattice points from $Q\T$ with relatively prime coordinates that generate $\kk$.

We conclude this section by mentioning that the symmetric role played
by numerators and denominators in \eqref{eqf} assures that just about
any statement, such as Theorems  1-4 below, holds identically with the word
'denominator' replaced by the word 'numerator'.

\section{The finite probability $\rho_Q^r(\cc,\d)$}

Let $\N_Q^r(\cc,\d)$ be the number of r-tuples of denominators of consecutive
fractions in $\FQ$ that are congruent with $\cc$ modulo $\d$, that is 
$\N_Q^r(\cc,\d)$ is the cardinality of the set
	\begin{equation*}
		\G_Q^r(\cc,\d)=\big\{(q_1,\dots,q_r)\in\cD^r_Q\colon\
		q_j\equiv c_j \pmod \d, \text{ for } j=1,\dots, r\big\}\,.
	\end{equation*}
Then we ought to to find the proportion
	\begin{equation*}
		\rho_Q^r(\cc,\d):=\frac{\N_Q^r(\cc,\d)}{\#\F_Q-(r-1)}		
	\end{equation*}
and probe the existence of a limit of $\rho_Q^r(\cc,\d)$, as $Q\to\infty$.

In order to find $\N_Q^r(\cc,\d)$, we need to estimate the number
of points with integer coordinates having a certain parity, and
belonging to different domains. This is the object of the next section.

\section{Lattice Points in Plane Domains}

Given a set $\Omega\subset\RR^2$ and  integers $0\le a,b<\d$, let
$N'_{a,b;\d}(\Omega)$ be the number of lattice points in $\Omega$ with
relatively prime coordinates congruent modulo $\d$ to $a$ and $b$,
respectively, that is, 
    \begin{equation*}
        N'_{a,b;\d}(\Omega)=\big\{
                        (m,n)\in\Omega\colon\
                                m\equiv a\pmod \d;\ n\equiv b\pmod \d;\
                                        \gcd(m,n)=1
                \big\}\,.
    \end{equation*}
Notice that $N'_{a,b;\d}=0$ if $\gcd(a,b)>1$, so  we may
assume that $a,b$ are relatively prime.
\begin{lemma}\label{LNO}
        Let $R>0$ and let $\Omega\subset\RR^2$ be a convex set of
diameter $\le R$. Let $\d$ be a positive integer and let $0\le a,b<\d$,
with $\gcd(a,b)=1$. Then
    \begin{equation}\label{eqNQ}
        N'_{a,b;\d}(\Omega)=\frac 6{\pi^2 \d^2}\prod_{p\mid \d}\left(1-\frac
                        1{p^2}\right)^{-1}
                        \Area(\Omega)+O\big(R\log R\big)\,.
    \end{equation}
\end{lemma}
\begin{proof}
As $\d$ is fixed and $R$ can be taken large enough, we may assume that
$\Omega\subset [0,R]^2$. Removing the coprimality condition by
M\"obius summation, we have:
    \begin{equation*}%
        \begin{split}
                N'_{a,b;\d}(\Omega)
                =&\sum_{\substack{(m,n)\in\Omega\\ m\equiv a\smod \d\\
                        n\equiv b \smod \d\\ \gcd(m,n)=1 }}1
                =\sum_{\substack{(m,n)\in\Omega\\ m\equiv a\smod \d\\
                        n\equiv b \smod \d }}
                        \sum_{\substack{e\mid m\\ e\mid n}}\mu(e)%
                =\sum_{e=1}^R\mu(e)
                \sum_{\substack{(m,n)\in\Omega\\ m\equiv a\smod \d\\
                        n\equiv b \smod \d\\ e\mid m,\ e\mid n }}1\,.
        \end{split}
    \end{equation*}
In the last sum, the last four conditions can be rewritten as:
$m=em_1$, $n=en_1$, $em_1\equiv a\pmod \d$, $en_1\equiv b\pmod \d$, for
some integers $m_1,n_1$. Since $\gcd(a,b)=1$, it follows that
$(e,\d)=1$. Let $e^{-1}$ be the inverse of $e \bmod \d$. Then,
    \begin{equation*}
        \begin{split}
                N'_{a,b;\d}(\Omega)
                        =\sum_{1\le e\le R}\mu(e)
                \sum_{\substack{(m_1,n_1)\in\frac 1e\Omega\\ m_1\equiv e^{-1}a\smod \d\\
                        n_1\equiv e^{-1}b \smod \d }}1\,.
        \end{split}
    \end{equation*}
Here the inner sum is equal to
$\frac 1{\d^2}\Area\big(\frac 1e \Omega\big)
        +O\Big(\length\Big(\partial\big(\frac 1e\Omega\big)\Big)\Big)$,
therefore
    \begin{equation*}
        \begin{split}
                N'_{a,b;\d}(\Omega)
                        =\frac 1{\d^2} \Area(\Omega)
                        \sum_{\substack{1\le e\le R\\
                                \gcd(e,\d)=1}}\frac{\mu(e)}{e^2}
                                +O\big(R\log R\big)\,.
        \end{split}
    \end{equation*}
Completing the last sum and the Euler product for the Dirichlet
series, we have
    \begin{equation*}
        \begin{split}
                \sum_{\substack{1\le e\le R\\
                        \gcd(e,\d)=1}}\frac{\mu(e)}{e^2}
                =\sum_{\substack{1\le e\le \infty\\
                        \gcd(e,\d)=1}}\frac{\mu(e)}{e^2}        +O(R^{-1})
                =\prod_p\left(1-\frac{1}{p^2}\right)
			\prod_{p\mid \d}\left(1-\frac{1}{p^2}\right)^{-1}
                        +O(R^{-1})\,.
        \end{split}
    \end{equation*}
This completes the proof of the lemma, since the first product is
equal to $6/\pi^2$.
\end{proof}

In particular, Lemma~\ref{LNO} provides an estimation of the cardinality of $\F_Q(c,d)$.

\begin{lemma}\label{LFcd}
        Let $d\ge 1$ and $0\le c\le d-1$ be integers. Then
    \begin{equation*}
        \#\F_Q(c,d)=\frac {3\nu(c,d)}{\pi^2 \d^2}\prod_{p\mid \d}
		\left(1-\frac 1{p^2}\right)^{-1} Q^2
                        +O\big(Q\log Q\big)\,,
    \end{equation*}
where $\nu(c,d)=\frac{\varphi(c)}{c}d+O(c)$ is the number of positive
integers $\le d$ that are relatively prime with $c$.
\end{lemma}
\begin{proof}
Let $\Omega=\T_Q$ be the triangle with vertices $(Q,0)$, $(Q,Q)$,
$(0,Q)$. Then, via what we know from Section~\ref{SecPere}, counting
only by the first component of a pair of consecutive elements of
$\FQ$, we find that  $\#\F_Q(c,d)$ is the number of lattice points
$(a,b)\in\T_Q$ with $\gcd(a,b)=1$, $a\equiv c\pmod d$ and no other
condition on $b$. Then the proof follows, since by Lemma~\ref{LNO}, we get
    \begin{equation*}
	\begin{split}
        \#\F_Q(c,d)&=\sum_{\substack{b=1\\ \gcd(b,d)=1}}^d N'_{c,b;\d}(\T_Q)\\
		&=\frac 6{\pi^2 \d^2}\prod_{p\mid \d}
		\left(1-\frac 1{p^2}\right)^{-1}
			\sum_{\substack{b=1\\ \gcd(b,d)=1}}^{d} \frac{Q^2}{2}
                        +O\big(Q\log Q\big)\,.	
	\end{split}
    \end{equation*}
\end{proof}
In particular, by Lemma~\ref{LFcd}, we find that 
$\#\FQ=\frac{3}{\pi^2}Q^2+O(Q\log Q)$ and the cardinality of the
subsets of fractions with odd and even denominators are 
$\#\FQo=\frac{2}{\pi^2}Q^2+O(Q\log Q)$ and
$\#\FQo=\frac{1}{\pi^2}Q^2+O(Q\log Q)$, respectively. 

\section{The densities $\rho^r_Q(\cc,\d)$ and $\rho^r(\cc,\d)$}\label{SectionPrinc}
We start with two particular cases. First we let $r=1$. Then
$\cc=c_1$ and, making use of Lemma~\ref{LFcd}, we obtain
	\begin{equation}\label{eqre1}
		\begin{split}
		\rho_Q^1(c_1,\d)&=\frac{\N_Q^1(c_1,\d)}{\#\F_Q}
			=\frac{\#\F_Q(c_1,\d)}{\#\F_Q}\\
		&=\frac {\nu(c_1,\d)}{\d^2}\prod_{p\mid \d}
		\left(1-\frac 1{p^2}\right)^{-1} +O\big(Q^{-1}\log Q\big)\,.
		\end{split}
	\end{equation}
Next let $r=2$. Now the numerator of $\rho^2_Q(\cc,\d)$ is the number of
lattice points from $\T_Q$ with relatively prime coordinates and
congruent with $\cc=(c_1,c_2)$, also. This is exactly the number counted in
Lemma~\ref{LNO}, thus we get: 
	\begin{equation}\label{eqre2}
		\begin{split}
		\rho_Q^2(c_1,c_2)&=\frac{\N_Q^2(c_1,c_2;\d)}{\#\F_Q-1}
			=\frac{N'_{c_1,c_2;\d}(\T_Q)}{\#\F_Q-1}\\
		&=\frac {1}{\d^2}\prod_{p\mid \d}
		\left(1-\frac 1{p^2}\right)^{-1} +O\big(Q^{-1}\log Q\big)\,.
		\end{split}
	\end{equation}

Now we assume that $r\ge3$ and $\cc=(c_1,\dots,c_r)$. We denote by
$\K^{r-2}(\cc,\d)$ the set of all 
$\kk=(k_1,\dots,k_{r-2})\in\A^{r-2}$ that correspond to $r$-tuples of 
consecutive denominators that are congruent to $\cc\pmod \d$, that is
  \begin{equation*}
     \K^{r-2}(\cc,\d)=\left\{\kk\in \A^{r-2}\colon\ 
        \begin{array}{l}  \exists Q\ge 1,\ \exists (q',q'')\in\T_Q,\ \gcd(q',q'') = 1,
                 \\ \displaystyle
	        \kk(q',q'')=\kk,\ \qq(q',q'')\equiv \cc\pmod\d
       \end{array} \right\}.
  \end{equation*}  
The natural continuum envelope of the generators of
$\kk\in\K^{r-2}(\cc,\d)$ is the set  
	\begin{equation*}
		\begin{split}
		\E_Q(\cc,\d)&=\big\{(x,y)\in\T_Q\colon\ 
		\kk(x,y)\in\K^{r-2}(\cc,\d)\big\}\\
			&=\bigcup_{\kk\in\K^{r-2}(\cc,\d)}Q\T[\kk]\,,
		\end{split}
	\end{equation*}
the union being disjoint.
Here are the lattice points that we have to count:
	\begin{equation*}
		\begin{split}
		\mathcal G_Q(\cc,\d)&=\left\{
					(q',q'')\in\T_Q\colon\ 
		        \begin{array}{l}  
		\gcd(q',q'')=1,\
		\kk(q'/Q,q''/Q)\in\K^{r-2}(\cc,\d),
	                 \\ \displaystyle
		q_1\equiv c_1\pmod\d,\ q_2\equiv c_2\pmod\d 
			       \end{array} \right\}\\
			&=\bigcup_{\kk\in\K^{r-2}(\cc,\d)}
			\big\{(q',q'')\in Q\T[\kk]\colon\ 
		\gcd(q',q'')=1,\ 
		(q_1,q_2)\equiv (c_1,c_2)\pmod\d\big\}\,.
		\end{split}
	\end{equation*}
Then
	\begin{equation*}
		\N^r_Q(\cc,\d)=N'_{c_1,c_2;\d}(\E_Q(\cc,\d))
		=\sum_{\kk\in\K^{r-2}(\cc,\d)}N'_{c_1,c_2;\d}\big(Q\T[\kk]\big)\,.
	\end{equation*}
Applying Lemma~\ref{LNO}, this can be written as
	\begin{equation}\label{eqF0}
		\begin{split}
		\N^r_Q(\cc,\d)=\frac{6Q^2}{\pi^2\d^2}
		\prod_{p\mid \d}\left(1-\frac 1{p^2}\right)^{-1} 
			&\sum_{\kk\in\K^{r-2}(\cc,\d)}\Area\big(\T[\kk]\big)\\
		&+O\Big(Q\log Q\sideset{}{'}\sum_{\kk\in\K^{r-2}(\cc,\d)}
				\length\big(\partial\T[\kk]\big)	
				\Big)\,.
		\end{split}
	\end{equation}
The prime attached to the sum in the error term indicates that in the
summation are excluded those vectors $\kk$ for which $\overline{Q\T[\kk]}$,
the adherence of $Q\T[\kk]$, contains no lattice points.  
We remark that the same estimate applies when $r=1,2$, as noticed in
the beginning of this section.
Depending on $r$, $\cc$ and $\d$, the nature of the sums in
\eqref{eqF0} may be different. In any case the first series is finite,
being always bounded from above by $\Area(\T)=1/2$. In the second
series, most of the times are summed infinitely many terms, but one
expects that its rate of convergence is small 
enough to assure a limit for $\rho^r_Q(\cc,\d)$ as $Q\to\infty$.
We summarize the result in the following theorem.

\begin{theorem}\label{Theorem1}
	Let $r\ge 3$, $\d\ge 2$ and $0\le c_1,\dots, c_r\le \d-1$ be
integers. Then
	\begin{equation}\label{eqF}
		\rho^r_Q(\cc,\d)=\frac{2}{\d^2}
		\prod_{p\mid \d}\left(1-\frac 1{p^2}\right)^{-1} 
			\sum_{\kk\in\K^{r-2}(\cc,\d)}\Area\big(\T[\kk]\big)
			+O\big(L_Q\big(\K^{r-2}(\cc,\d)\big)Q^{-1}\log Q\big)\,,
	\end{equation}
where $L_Q\big(\K^{r-2}(\cc,\d)\big)$ is the sum of the perimeters of all
polygons $Q\T[\kk]$ with $\kk\in\K^{r-2}(\cc,\d)$, and having the property
that $\overline{Q\T[\kk]}$ contains lattice points.
\end{theorem}

In order to have a limit of the ratio $\rho^r_Q(\cc,\d)$,
we need at least to know that the limit
$L_Q\big(\K^{r-2}(\cc,\d)\big)/ Q^{2}\to 0$, as 
$Q\to\infty$, exists. A few initial checks leads one to expect a much
stronger estimate to be true. Indeed, the authors~\cite{Recursive}
have proved that for any fixed $r\ge 1$, the number of $\kk\in\A^r$ with all
components $\le Q$ is $rQ+O(r)$. Moreover, for $Q$ sufficiently large,
the polygons corresponding to $\kk$ with  components larger than $Q$
(and only one  components can be so\,!) form at most two connected
regions in $\T$. These are smaller 
and smaller when $Q\to\infty$ and tend to one limit point, $(1,0)$,
in the case $r=1$,  and to two limit points, $(1,1)$ and $(1,0)$, when
$r\ge 2$.  
Regarding the perimeters of $\T[\kk]$, in \cite{Recursive} and
\cite{Tessellations}  we have shown that
        \begin{equation}\label{eqPerimetrele}
		\sum_{\substack{\kk\in\A^r\\ 1\le k_1,\dots,k_r\le Q}}
		\length \big(\partial \T[\kk]\big)
		\ll r\log Q \,,
	\end{equation}
and consequently this is also an upper bound for
$L_Q\big(\K^{r-2}(\cc,\d)\big)$. One can find in ~\cite{Tessellations}
more details on the tessellation of $\T$ formed by the polygons
$\T[\kk]$ with $\kk$ of the same order $r$. For instance, the polygons
$\T[\kk]$, whose vectors are excepted in the domain of summation in
\eqref{eqPerimetrele} have the  form  
$\kk=(\underbrace{2,\dots,2,1}_{s\ \text{components}},k,
		\underbrace{1,2,\dots,2}_{t\ \text{components}})$,
with $s,t\ge 0$,  $s+1+t=r$ and $k> Q$, are quadrangles whose vertices
are given by formulas that are calculated explicitly. These are
exactly the polygons that have the main contribution (by their number) in the
estimate~\eqref{eqPerimetrele}, since the number of the remaining ones
is always finite. Employing these information in
Theorem~\ref{Theorem1}, we obtain our main result.  
\begin{theorem}\label{Theorem2}
	Let $r\ge 3$, $\d\ge 2$ and $0\le c_1,\dots, c_r\le \d-1$ be
integers. Then
	\begin{equation}\label{eqF2}
		\begin{split}
		\rho^r_Q(\cc,\d)&=\frac{2}{\d^2}
		\prod_{p\mid \d}\left(1-\frac 1{p^2}\right)^{-1} 
			\sum_{\kk\in\K^{r-2}(\cc,\d)}\Area\big(\T[\kk]\big)
				+O\big(rQ^{-1}\log^2 Q\big)\,.
		\end{split}
	\end{equation}
\end{theorem}

\begin{corollary}\label{Corollary1}
	Let $r\ge 3$, $\d\ge 2$ and $0\le c_1,\dots, c_r\le \d-1$ be
integers. Then, there exists the limit
$\rho^r(\cc,\d):=\lim_{Q\to\infty}\rho^r_Q(\cc,\d)$, and  
	\begin{equation*}
		\begin{split}
		\rho^r(\cc,\d)&=\frac{2}{\d^2}
		\prod_{p\mid \d}\left(1-\frac 1{p^2}\right)^{-1} 
			\sum_{\kk\in\K^{r-2}(\cc,\d)}\Area\big(\T[\kk]\big)\,.
		\end{split}
	\end{equation*}
\end{corollary}

\section{On short intervals}\label{secShortIntervals}
We now turn to see what changes occur when we treat the same problem
on an arbitrary interval. By the periodicity modulo an interval of
length $1$ of consecutive denominators of fractions in $\F^\I_Q$, we
may reduce to consider only shorter intervals. Thus, in the following
we assume that the interval $\I\subseteq [0,1]$, of positive length, is fixed. 
In the notations introduced above for different sets, we will use an
additional superscript $\I$ with the significance that their elements
correspond to fractions from $\I$.

By the fundamental relation \eqref{eqf}, we find that if
$(a'/q',a''/q'')\in\cF_Q^{2}$ then
	\begin{equation}\label{eqCondInv}
		a''/q''\in\I \iff \overline{q'}\in q''\I\,,
	\end{equation}
in which $\overline{q'}\in[0,q'']$ is the inverse of $q'\pmod {q''}$.
Let
	\begin{equation*}
		\G_Q^{\I,r}(\cc,\d)=\big\{(q_1,\dots,q_r)\in\cD^{\I,r}_Q\colon\
		q_j\equiv c_j \pmod \d, \text{ for } j=1,\dots, r\big\}\,
	\end{equation*}
and $\N_Q^{\I,r}(\cc,\d)=\#\G_Q^{\I,r}(\cc,\d)$. 
Then, our task is to estimate
	\begin{equation}\label{eqrhoI}
          \rho_Q^{\I,r}(\cc,\d)
		:=\frac{\N_Q^{\I,r}(\cc,\d)}{\#\F_Q^{\I}-(r-1)}\,.	
	\end{equation}
Now we consider the set
	\begin{equation*}
		\begin{split}
		\mathcal G_Q^\I(\cc,\d)&=
			\left\{(q',q'')\in\T_Q\colon\ 
			        \begin{array}{l}  
			\overline{q'}\in q''\I,\
				(q',q'')\equiv (c_1,c_2) \pmod {\d},\
                 \\ \displaystyle
			\gcd(q',q'')=1,\ \kk(q',q'')\in\K^{r-2}(\cc,\d)
       \end{array} \right\}\\
			&=\bigcup_{\kk\in\K^{r-2}(\cc,\d)}
				\left\{(q',q'')\in Q\T[\kk]\colon\ 
			        \begin{array}{l}  
			\overline{q'}\in q''\I,\ \gcd(q',q'')=1,
                 \\ \displaystyle
			(q',q'')\equiv (c_1,c_2) \pmod {\d}\
			       \end{array} \right\}.
		\end{split}
	\end{equation*}
Then
	\begin{equation}\label{eqI3}
		\N^{\I,r}_Q(\cc,\d)=\#\mathcal G_Q^\I(\cc,\d)+O(1)\,.
	\end{equation}

For a given plane domain $\Omega$, we use the following notation
  \begin{equation*}
     N^{',\I}_{c_1,c_2;\d}(\Omega):=\#\left\{(q',q'')\in\Omega\cap\NN^2\colon\ 
        \begin{array}{l} \overline{q'}\bmod q''\in q''\I,\ \gcd(q',q'') = 1,
                 \\ \displaystyle
                         (q_1,q_2)\equiv (c_1,c_2)\pmod\d
        \end{array} \right\}.
  \end{equation*}  
Then, \eqref{eqI3} becomes
	\begin{equation}\label{eqI4}
		\N^{\I,r}_Q(\cc,\d)=
	\sum_{\kk\in\K^{r-2}(\cc,\d)}N^{',\I}_{c_1,c_2;\d}\big(\T[\kk]\big)
			+O(1)\,.
	\end{equation}

The next lemma shows that $N^{',\I}_{c_1,c_2;\d}(\Omega)$ is, roughly,
$N^{'}_{c_1,c_2;\d}(\Omega)$ times the length of the interval $\I$.

\begin{lemma}\label{LemmaNc1c2I}
Let $Q>0$ and let $\Omega$ be a convex domain included in the triangle of
vertices $(Q,0); (Q,Q); (0,Q)$.
Let $\d$ be a positive integer and let $0\le a,b<\d$,
with $\gcd(a,b)=1$. 
Then, we have
   \begin{equation*}
	N^{',\I}_{a,b;\d}(\Omega)= |\I|\cdot N^{'}_{a,b;\d}(\Omega)
	 + O(Q^{3/2+\varepsilon})\,.
   \end{equation*}
\end{lemma}
\begin{proof}
We count the good points in $\Omega$ situated on horizontal lines with integer
coordinates. Thus, we have
  \begin{equation}\label{eqNFeI1}
        \begin{split}
	N^{',\I}_{c_1,c_2;\d}(\Omega)
		=\sum_{\substack{1\le q\le Q\\ q\equiv b\pmod\d}}\#\left\{
		x\in \Omega\cap\{ y=q\}\colon \ 
        \begin{array}{l} 
		\overline{x}\bmod q\in q\I,\ \gcd(x,q)=1,
                 \\ \displaystyle
		 x\equiv a\pmod\d \end{array} \right\}.
        \end{split}
  \end{equation}  
Employing exponential sums, the terms in the sum are equal to
  \begin{equation}\label{eqNFeI2}
        \begin{split}
        \sum_{\substack{x\in\Omega\cap\{ y=q\}\\ \gcd(x,q)=1\\
			x\equiv a \pmod \d}}\ \
        \sum_{y\in q\I}
        \frac 1q\sum_{k=1}^qe\Big(k\frac{y-\overline x}{q}\Big)
                &=\frac 1q\sum_{k=1}^q        
			\sum_{y\in q\I}e\Big(k\frac{y}{q}\Big)
	\sum_{\substack{x\in\Omega\cap\{ y=q\}\\ \gcd(x,q)=1\\
			x\equiv a \pmod \d}}\ 
	e\Big(k\frac{-\overline x}{q}\Big)\,.
        \end{split}
  \end{equation}  
We separate the terms in \eqref{eqNFeI2} in two groups. The first one
contains the terms with $k=q$ and the second is formed by all the others. 
The contribution of the terms from the first group will give the main
term in \eqref{eqNFeI1}, since 
  \begin{equation}\label{eqL1}
        \begin{split}
	\sum_{\substack{1\le q\le Q\\ q\equiv b\pmod\d}}
        \sum_{\substack{x\in\Omega\cap\{ y=q\}\\ \gcd(x,q)=1\\
			x\equiv a \pmod \d}}\ \
        \sum_{y\in q\I}        \frac 1q
	=|\I|\cdot N^{'}_{a,b;\d}(\Omega)+O(Q)\,.
        \end{split}
  \end{equation}  
It remains to estimate the size of the terms from the second group.
The second sum from the right-hand side of \eqref{eqNFeI2} is a
geometric progression that is bounded sharply by $\ll
\min(q|\I|,\norm{k/q}^{-1})$, where $\norm{\cdot}$ is  the distance to
the closest integer. The most inner sum is a Kloosterman-type sum. It
is incomplete, both on the length of the interval and on the
$x$-domain--an arithmetic progression. A standard procedure, using the
classic bound of Esterman~\cite{Esterman} and Weil~\cite{Weil}, gives
  \begin{equation*}
	\bigg|\sum_{\substack{x\in\Omega\cap\{ y=q\}\\ \gcd(x,q)=1\\
			x\equiv a \pmod \d}}\ 
	e\Big(k\frac{-\overline x}{q}\Big)\bigg|
	\le \sigma_0(q) (k,q)^{1/2}q^{1/2}(2+\log q)\,,
  \end{equation*}  
where $\sigma_l(q)$ is the sum of the $l$-th power of divisors of $q$.
Thus, the contribution to \eqref{eqNFeI1} of the terms from the second
group is
  \begin{equation}\label{eqL21}
        \begin{split}
	\ll
	\sum_{\substack{1\le q\le Q\\ q\equiv b\pmod\d}}
		\sigma_0(q) q^{1/2}(2+\log q)
		\frac 1q\sum_{k=1}^{q-1} (k,q)^{1/2} \min(q|\I|,\norm{k/q}^{-1})\,.    
        \end{split}
  \end{equation}  
Here the sum over $k$ is
  \begin{equation*}
        \begin{split}
		\sum_{k=1}^q   (k,q)^{1/2} \norm{k/q}^{-1}
		&=\sum_{g\mid q}\sum_{\substack{k=1 \\ (k,q)=g}}^q
			\frac{g^{1/2}}{\norm{k/q}}\\
	&=\sum_{g\mid q}g^{1/2}\sum_{k=1}^{\big[\frac{q-1}{2g}\big]}
		\frac{2q}{gk}	\\
	&\le 2 \sigma_{-1/2}(q)q(2+\log q)\,.
        \end{split}
  \end{equation*}  
On inserting this estimate in \eqref{eqL21} and using the fact that
$\sigma_l(q)\ll q^{\varepsilon}$, we see that the contribution of
terms from the second group is $\ll Q^{3/2+\varepsilon}$.
This completes the proof of the lemma.
\end{proof}

In particular, Lemma~\ref{LemmaNc1c2I} may be used to count the
fractions from an interval. Thus, we have
	\begin{equation}\label{eqFI}
		\#\F_Q^{\I}=|\I|\cdot\#\FQ+O\big(Q^{3/2+\varepsilon}\big)\,.
	\end{equation}

Finally, using the estimate given by Lemma~\ref{LemmaNc1c2I} in
\eqref{eqI4} and combining the result and \eqref{eqFI} in
\eqref{eqrhoI}, we obtain the following theorem.

\begin{theorem}\label{TheoremI}
	Let $r\ge 1$, $\d\ge 2$ and $0\le c_1,\dots, c_r\le \d-1$ be
integers. Then
	\begin{equation*}
		\rho^{\I,r}_Q(\cc,\d)= \rho^{r}_Q(\cc,\d)
			+O\big(Q^{-1/2}\log^2 Q\big)\,.
	\end{equation*}
\end{theorem}

As a consequence, it follows that independent on the interval, the limit
$\rho^{\I,r}(\cc,\d):=\lim_{Q\to\infty}\rho^{\I,r}_Q(\cc,\d)$ exists.

\begin{corollary}\label{CorollaryI}
	Let $r\ge 1$, $\d\ge 2$ and $0\le c_1,\dots, c_r\le \d-1$ be
integers. Then, for any interval $\I$ of positive length, the sequence
$\{\rho^{\I,r}_Q(\cc,\d)\}_{Q}$ has a limit $\rho^{\I,r}(\cc,\d)$, as
$Q\to\infty$, and 
	\begin{equation*}
		\rho^{\I,r}(\cc,\d)= \rho^{r}(\cc,\d)\,.
	\end{equation*}
\end{corollary}
\section{A few special cases}
We begin with the case $\d=2$. Then $\cc\in\{0,1\}^r$, that is we are looking
to the probability that 
$r$-tuples of consecutive denominators are odd or even in a prescribed
order. We may always suppose that there are no neighbor even
components of $\cc$, since in that case $\rho^r(\cc,2)=0$.

When $r=1$ or $r=2$, we already know from the beginning of
Section~\ref{SectionPrinc}, relations  \eqref{eqre1} and \eqref{eqre2},
the precise values of $\rho^r(\cc,2)$, while for 
$r\ge 3$, we get $\rho^r(\cc,2)$ from Theorem~\ref{Theorem1}.  

\begin{theorem}\label{TheoremDoi}
	Let $r\ge 1$ and $c_1,\dots, c_r\in\{0,1\}$. Then, there
exists the limit
$\rho^r(\cc,2)=\lim_{Q\to\infty}\rho^r_Q(\cc,2)$. Furthermore, we
have:
$\rho^1(0;2)=1/3$,
$\rho^1(1;2)=2/3$, $\rho^1(0,1;2)=\rho^1(1,0;2)=\rho^1(1,1;2)=1/3$, and
	\begin{equation}\label{eqThDoi}
		\begin{split}
		\rho^r(\cc,2)&=%
		\frac 23
		\sum_{\kk\in\K^{r-2}(\cc,2)}\Area\big(\T[\kk]\big), \quad
				\text{ for $r\ge 3$}\,.
		\end{split}
	\end{equation}
\end{theorem}

We have calculated the sums from the right-hand side of
\eqref{eqThDoi} in a few cases. Here they are. 
First we remark that when $\kk$ has only one component (i.e. it
corresponds to $3$-tuples of consecutive denominators), the areas are:
$\Area\big(\T[1]\big)=1/6$ and
$\Area\big(\T[k]\big)=4/\big(k(k+1)(k+2)\big)$, for 
$k\ge 2$. Then, employing the sum of the Leibniz series, we obtain: 
   \begin{equation*}
	\begin{split}
		\rho^{3}(1,1,1;2)&= \frac 23 
		\sum_{\substack{k\ge 1\\  k \text{ even}}}\Area\big(\T[k]\big)
		=2-\frac 83 \log 2\approx 0.15160\,;\\
		\rho^{3}(1,1,0;2)&= \frac 23 
		\sum_{\substack{k\ge 1\\  k \text{ odd}}}\Area\big(\T[k]\big)
		=\frac 83 \log 2-\frac 53\approx 0.18172\,;\\
		\rho^{3}(1,0,1;2)&= \frac 23 
		\sum_{\substack{k\ge 1}}\Area\big(\T[k]\big)
		=\frac 13\,\approx 0.33333;\\
		\rho^{3}(0,1,1;2)&= \rho^{3}(1,1,0;2)
			=\frac 83 \log 2-\frac 53\approx 0.18172\,;\\
		\rho^{3}(0,1,0;2)&= \frac 23 
		\sum_{\substack{k\ge 1\\  k \text{ even}}}\Area\big(\T[k]\big)
		=2-\frac 83 \log 2\approx 0.15160\,.\\
	\end{split}
   \end{equation*}
Thus, out of the $8$ possible vectors $\cc$, only $5$ are not trivial
(since the others have two neighbor even denominators, so
$\rho^{3}(0,0,1;2)=\rho^{3}(1,0,0;2)=\rho^{3}(0,0,0;2)=0$). Furthermore two of
them form a couple with the same occurring probability, because their
components are mirrorly reflected of one another.

For longer sequences $\cc$, the sum from the right-hand side of
\eqref{eqThDoi} involves the Leibniz series  again, more precisely its smaller
 and smaller remainder. This happens because more and more
$\kk$'s with all components small belong to $\K^{r-2}(\cc,2)$, while
those with at least one component $k$, say, passing over a certain
magnitude have the property that
$\Area\big(\T[\kk]\big)=\Area\big(\T[k]\big)$. We 
mention here only that
   \begin{equation*}
	\begin{split}
		\Area\big(\T[1,k]\big)&=\Area\big(\T[k,1]\big)
		=\Area\big(\T[k]\big)\,,
			\quad\text{ for $k\ge 5$;}\\
		\Area\big(\T[1,k,1]\big)&=\Area\big(\T[k]\big)\,,
			\quad\text{ for $k\ge 5$;}\\
		\Area\big(\T[2,1,k]\big)&=\Area\big(\T[k,1,2]\big)
		=\Area\big(\T[k]\big)\,,
			\quad\text{ for $k\ge 9$.}
	\end{split}
   \end{equation*}

In the case $r=4$, there are $8$ nontrivial vectors $\cc$, of which
$5$ are essentially distinct (non mirror reflected of another). Here
are the probabilities with which they come about:
   \begin{equation*}
	\begin{split}
		\rho^{4}(1,1,1,1;2)&= \frac 23 
		\sum_{\substack{k,l\ge 1\\  k\text{ even, },l \text{ even}}}
			\Area\big(\T[k,l]\big)
		=\frac{23}{315}\approx 0.07301 \,;\\
		\rho^{4}(1,1,1,0;2)&= \frac 23 
		\sum_{\substack{k,l\ge 1\\  k \text{ even, }l \text{ odd}}}
			\Area\big(\T[k,l]\big)
		=\frac{607}{315}- \frac83\log 2\approx 0.07859\,;\\
		\rho^{4}(1,1,0,1;2)&= \frac 23 
		\sum_{\substack{k,l\ge 1\\  k \text{ odd }}}\Area\big(\T[k,l]\big)
		=\frac 83\log 2-\frac 53\,\approx 0.18172\,;\\
		\rho^{4}(1,0,1,1;2)&= \rho^{4}(1,1,0,1;2)
			=\frac 83\log 2-\frac 53\approx 0.18172\,;\\
		\rho^{4}(1,0,1,0;2)&= \frac 23 
		\sum_{\substack{k,l\ge 1\\  l \text{ even}}}\Area\big(\T[k,l]\big)
		=2-\frac 83 \log 2\approx 0.15160\,;\\
		\rho^{4}(0,1,1,1;2)&= \rho^{4}(1,1,1,0;2)
		=\frac{607}{315}- \frac83\log 2\approx 0.07859\,;\\
		\rho^{4}(0,1,1,0;2)&= \frac 23 
		\sum_{\substack{k,l\ge 1\\  k \text{ odd, }l \text{ odd}}}
			\Area\big(\T[k,l]\big)
		=\frac{16}{3}\log 2- \frac{1132}{315}\approx 0.10313\,;\\
		\rho^{4}(0,1,0,1;2)&= \rho^{4}(1,0,1,0;2)
		=2-\frac 83 \log 2\approx 0.15160\,.\\
	\end{split}
   \end{equation*}

When $r=5$, out of the $32$ vectors $\cc\in\{0,1\}^5$, only $13$ have
no neighbor even-even components and $9$ are essentially distinct. We
present the probabilities $\rho^5(\cc,2)$ in  Table~\ref{Table1}.

\smallskip
\setlength{\doublerulesep}{1pt}
\setlongtables
\begin{longtable}{|C|C|C|C|}
\caption{The probabilities $\rho^5(\cc,2)$. In the $\kk$-column, the
notations $e$, $o$ and $\forall$ mean that the sum from the right-hand
side of \eqref{eqThDoi} runs over all vectors with the correspondent
component even, odd, or whatever, respectively.}\label{Table1}\\ \hline
{\normalfont \cc} & {\normalfont \kk} & \rho^5(\cc,2)  &\mathrm{approximation} \\ \hhline{|====|}%
\endfirsthead
\multicolumn{3}{l}{\small\sl continued from previous page}\\ \hline
{\normalfont \cc} & {\normalfont \kk} & \rho^5(\cc,2)  &\mathrm{approximation} \\ \hhline{|====|}%
\endhead
\hline
\multicolumn{4}{r}{\small\sl continued on next page} \\ %
\endfoot
\hline
\endlastfoot
(1,1,1,1,1) & (e,e,e) & \frac{1}{21} & 0.04761 \\ \hline
(1,1,1,1,0) & (e,e,o) & \frac{8}{315} & 0.02539 \\ \hline
(1,1,1,0,1) & (e,o,\forall) & \frac{607}{315}-\frac{8}{3}\log 2 &  0.07859\\ \hline
(1,1,0,1,1) & (o,\forall,o) & \frac{4441}{38610} &  0.11502\\ \hline
(1,1,0,1,0) & (o,\forall,e) & \frac{8}{3}\log 2-\frac{68791}{38610} & 0.06670 \\ \hline
(1,0,1,1,1) & (\forall,o,e) & \frac{607}{315}-\frac{8}{3}\log 2 &  0.07859\\ \hline
(1,0,1,1,0) & (\forall,o,o) & \frac{16}{3}\log 2 -\frac{1132}{315}&  0.10313\\ \hline
(1,0,1,0,1) & (\forall,e,\forall) & 2-\frac{8}{3}\log 2 &  0.15160\\ \hline
(0,1,1,1,1) & (o,e,e) & \frac{8}{315} & 0.02539\\ \hline
(0,1,1,1,0) & (o,e,o) & \frac{599}{315}-\frac{8}{3}\log 2&  0.05319\\ \hline
(0,1,1,0,1) & (o,o,\forall) & \frac{16}{3}\log 2-\frac{1132}{315} &  0.10313\\ \hline
(0,1,0,1,1) & (e,\forall,o) & \frac{8}{3}\log 2-\frac{68791}{38610} &  0.06670\\ \hline
(0,1,0,1,0) & (e,\forall,e) & \frac{146011}{38610}-\frac{16}{3}\log 2&  0.08490\\ \hline
\end{longtable}
\normalsize

Many patterns of consecutive denominators extend without bound as $Q$
gets large.  We mention here the one with all components equal modulo
$\d$. Let $\d\ge 2$ and $0\le c\le \d -1$. The condition of neighborship
produces the constrain $\gcd(c,\d)=1$. Let
$\cc=(c,\dots,c)$ be the vector with all the $r$ components equal to
$c$. Remarkably, when $r\ge 5$,  there exists only one $\kk$ which
accommodates the appearing in $\FQ$ of sequences of consecutive
denominators that are congruent with $\cc$ modulo $\d$. This is the
vector $\kk=(2,\dots,2)$ with  $r-2$ components all equal with $2$.
The corresponding polygon is the quadrangle (below we refer to formulas
proved in~\cite{Tessellations})
   \begin{equation*}
	\T_{r-2}[2,\dots,2]
	=\Big\{(1,1);\ \Big(\frac{r-2}{2r-3},\frac{r-1}{2r-3}\Big);\
		\Big(\frac 12,\frac 12\Big);\ \Big(1, \frac{2r-4}{2r-3}\Big)
	\Big\}\,,\quad\text{for $r\ge 3$}
   \end{equation*}
and its area is
   \begin{equation*}
	\Area\big(\T_{r-2}[2,\dots,2]\big)
	=\frac{1}{2(2r-3)}\,,\quad\text{for $r\ge 3$.}   
   \end{equation*}
(Since $r$, the number of components of $\kk$ is essential in the
formulae, in order to indicate precisely its size, we write
$\T_r[\kk]$ instead of $\T[\kk]$.)  
Then, Corollary~\ref{Corollary1} yields the following result.
\begin{corollary}\label{Corollarycconst}
	Let $r\ge 5$, $\d\ge 2$, and let $0\le c\le \d-1$ with $\gcd(c,\d)=1$.
Then
	\begin{equation*}%
		\rho^r(c,\dots,c;\d)=\frac{1}{\d^2(2r-3)}
		\prod_{p\mid \d}\left(1-\frac 1{p^2}\right)^{-1} \,.
	\end{equation*}
\end{corollary}
In particular, Corollary~\ref{Corollarycconst} gives the probability
to find $r\ge 5$ odd consecutive denominators:
   \begin{equation}\label{eqc11}
	\rho^r(\underbrace{1,\dots,1}_{r\text{ ones}};2)
	=\frac{1}{3(2r-3)}\,,\quad\text{for $r\ge 5$}\,.
   \end{equation}

The same pattern boarded on either side by an even denominator is very
similar. Indeed, if  
$\cc=(0,\underbrace{1,\dots,1}_{r-1 \text{ ones}})$ and $r\ge 6$,
there exists only two corresponding vectors:
$\kk=(1,\underbrace{2,\dots,2}_{r-3 \text{ twos}})$ and  
$\kk=(3,\underbrace{2,\dots,2}_{r-3 \text{ twos}})$,
while the mirror reflected case 
$\cc=(\underbrace{1,\dots,1}_{r-1 \text{ ones}},0)$ corresponds to
$\kk=(\underbrace{2,\dots,2}_{r-3 \text{ twos}},1)$ and
$\kk=(\underbrace{2,\dots,2}_{r-3 \text{ twos}},3)$. 
In all four cases $\T[\kk]$ is a triangle:
   \begin{equation*}
	\begin{split}
	\T[1,\underbrace{2,\dots,2}_{r-3 \text{ twos}}]
	&=\Big\{(0,1);\ \Big(\frac{1}{2r-3},\frac{2r-4}{2r-3}\Big);\
		 \Big(\frac{1}{2r-5}, 1\Big)
	\Big\},\quad\text{for $r\ge 3$}\,;\\
	\T[3,\underbrace{2,\dots,2}_{r-3 \text{ twos}}]
	&=\Big\{\Big(1,\frac 12\Big);\ \Big(1,\frac{r-1}{2r-3}\Big);\
		 \Big(\frac{2r-7}{2r-5}, \frac{r-3}{2r-5}\Big)
	\Big\},\quad\text{for $r\ge 5$}\,;\\
	\T[\underbrace{2,\dots,2}_{r-3 \text{ twos}},1]
	&=\Big\{(1,1);\ \Big(\frac{r-3}{2r-5},\frac{r-2}{2r-5}\Big);\
		 \Big(\frac{r-2}{2r-3}, \frac{r-1}{2r-3}\Big)
	\Big\},\quad\text{for $r\ge 3$}\,;\\
	\T[\underbrace{2,\dots,2}_{r-3 \text{ twos}},3]
	&=\Big\{\Big(\frac 12,\frac 12\Big);\ \Big(1,\frac{2r-6}{2r-5}\Big);\
		 \Big(1, \frac{2r-4}{2r-3}\Big)
	\Big\},\quad\text{for $r\ge 5$}\,;\\
	\end{split}
   \end{equation*}
and they have the same area:
   \begin{equation*}
	\begin{split}
	\Area\big(\T[1,\underbrace{2,\dots,2}_{r-3 \text{ twos}}]\big)
	&=\Area\big(\T[\underbrace{2,\dots,2}_{r-3 \text{ twos}},1]\big)\\
        &=\Area\big(\T[3,\underbrace{2,\dots,2}_{r-3 \text{ twos}}]\big)
	=\Area\big(\T[\underbrace{2,\dots,2}_{r-3 \text{ twos}},3]\big)\\
	&=\frac{1}{2(2r-5)(2r-3)}\,,\quad\text{for $r\ge 5$}\,.\\   
	\end{split}
   \end{equation*}
Then, by Corollary~\ref{Corollary1}, we get
   \begin{equation*}
	\rho^r(0,\underbrace{1,\dots,1}_{r-1\text{ ones}};2)
	=\rho^r(\underbrace{1,\dots,1}_{r-1\text{ ones}},0;2)
	=\frac{2}{3(2r-5)(2r-3)}\,,\quad\text{for $r\ge 6$}\,.   
   \end{equation*}

We conclude with an analogue example on the side of $\kk$'s. 
The question we address is whether beside $\kk=(2,\dots,2)$, there
exists another $\kk$ with all components equal, which extends without
bound. This happens, but only modulo some $\d\ge 2$, namely
$\d=3$, and $\kk$ being a series of ones intercalated by fours. As
$\kk$ and its pal, the one with components in reversed order, satisfy the
demands at the same time, the components of these vectors depends on
the parity of $r$. Precisely, for any $r\ge 1$, they are:
	\begin{equation*}%
		\begin{split}
			\T_{2r}[1,4,\dots,1,4]&=\Big\{
		\Big(\frac{1}{3},\frac{2}{3}\Big)\,; 
		\Big(\frac{3r}{6r-1},1\Big)\,;
		\Big(\frac{1}{2},1\Big)\,; 
		\Big(\frac{2r}{6r+1},\frac{4r+1}{6r+1}\Big)
				\Big\}\,;\\
			\T_{2r+1}[1,4,\dots,4,1]&=\Big\{
		\Big(\frac 13,\frac{2}{3}\Big)\,; 
		\Big(\frac{3r+1}{6r+1},1\Big)\,;
		\Big(\frac{1}{2},1\Big)\,; 
		\Big(\frac{2r}{6r+1},\frac{4r+1}{6r+1}\Big)
				\Big\}\,;\\
			\T_{2r}[4,1\dots,4,1]&=\Big\{
		\Big(1,\frac{1}{2}\Big)\,; 
		\Big(\frac{4r-1}{6r-1},\frac{2r}{6r-1}\Big)\,;
		\Big(\frac{2}{3},\frac{1}{3}\Big)\,; 
		\Big(1,\frac{3r}{6r+1}\Big)
				\Big\}\,;\\
			\T_{2r+1}[4,1,\dots,1,4]&=\Big\{
		\Big(1,\frac{1}{2}\Big)\,; 
		\Big(\frac{4r+3}{6r+5},\frac{2r+2}{6r+5}\Big)\,;
		\Big(\frac{2}{3},\frac{1}{3}\Big)\,; 
		\Big(1,\frac{3r+2}{6r+5}\Big)
				\Big\}\,.
		\end{split}
	\end{equation*}
And here are their areas:
	\begin{equation}\label{eqAreas14}
		\begin{split}
	\Area\big(\T_{2r}[1,4,\dots,1,4]\big)&=\frac{r}{36r^2-1}\,;\\
	\Area\big(\T_{2r+1}[1,4,\dots,4,1]\big)&=\frac{1}{36r+6}\,;\\
	\Area\big(\T_{2r}[4,1\dots,4,1]\big)&=\frac{r}{36r^2-1}\,;\\
	\Area\big(\T_{2r+1}[4,1,\dots,1,4]\big)&=\frac{1}{36r+30}\,.
		\end{split}
	\end{equation}
One of the patterns suited by this $\kk$'s is formed by sequences of
denominators that are congruent modulo $3$ with a series of repeated
ones and twos. For these, we obtain the following probabilities.
\begin{corollary}\label{CorollaryUnuDoi}
	For any $r\ge 4$, we have:
	\begin{equation*}
		\begin{split}
		\rho^{2r}(1,2,\dots,1,2;3)&=\rho^{2r}(1,2,\dots,1,2;3)=
	\frac{r-1}{72(r-1)^2-2},\\
		\rho^{2r+1}(1,2,\dots,2,1;3)&=\rho^{2r+1}(2,1,\dots,1,2;3)=
	\frac{9r+4}{8(9r+1)(9r+7)}\,.\\
		\end{split}
	\end{equation*}
\end{corollary}
\begin{proof}
By Corollary~\ref{Corollary1}, we have:
	\begin{equation*}
		\begin{split}
		\rho^{2r}(1,2,\dots,1,2;3)&=\frac 14\Big(
	\Area\big(\T_{2(r-1)}[1,4,\dots,1,4]\big)+
		\Area\big(\T_{2(r-1)}[4,1,\dots,4,1]\big)\Big)\,,\\
		\rho^{2r+1}(1,2,\dots,2,1;3)&=\frac 14\Big(
	\Area\big(\T_{2r-1}[1,4,\dots,4,1]\big)+
		\Area\big(\T_{2r-1}[4,1,\dots,1,4]\big)\Big)\,,\\
		\end{split}
	\end{equation*}
and the same relations apply for $\rho^{2r}(2,1,\dots,2,1;3)$ and
$\rho^{2r+1}(2,1,\dots,1,2;3)$, respectively. Then the corollary
follows using the formulae from \eqref{eqAreas14}.
\end{proof}

\end{document}